\documentclass{gtart}

\def\ifplaintex{\expandafter\ifx\csname documentclass\endcsname\relax}

\def\gtp{{\mathsurround=0pt\it $\cal G\mskip-2mu$eometry \&\ 
$\cal T\!\!$opology $\cal P\!$ublications}}  

\def\recd{{\small Received:\qua\receiveddate\ifx\reviseddate\relax
\else\qquad Revised:\qua\reviseddate\fi\par}} 

\def\lognumber#1{\def\thelognumber{#1}}
\def\volumenumber#1{\def\thevolumenumber{#1}}
\def\volumeyear#1{\def\thevolumeyear{#1}}
\def\papernumber#1{\def\thepapernumber{#1}}
\def\pagenumbers#1#2{\def\startpage{#1}\def\finishpage{#2}}
\def\published#1{\def\publishdate{#1}}

\def\received#1{\def\receiveddate{#1}}

\def\accepted#1{\def\accepteddate{#1}}

\def\asciiaddress#1{\def\theasciiaddress{#1}}
\def\asciiemail#1{\def\theasciiemail{#1}}

\long\def\asciiabstract#1{\long\def\theasciiabstract{#1}}
\def\asciikeywords#1{\def\theasciikeywords{#1}}

\let\\\par\let\thelognumber\relax\let\thevolumenumber\relax
\let\thepapernumber\relax\let\thevolumeyear\relax\let\startpage\relax
\let\finishpage\relax\let\publishdate\relax\let\receiveddate\relax
\let\reviseddate\relax\let\accepteddate\relax\let\theasciititle\relax
\let\theasciiauthors\relax\let\theasciiaddress\relax
\let\theasciiabstract\relax\let\theasciikeywords\relax

\let\theasciiemail\relax

\ifplaintex
\font\logobig=cmssbx10 scaled 3836
\font\logomed=cmssbx10 scaled 2557
\else
\font\logobig=cmssbx10 scaled 4200
\font\logomed=cmssbx10 scaled 2800
\fi

\long\def\makeagttitle{   
\count0=\startpage
\agt\hfill      
\hbox to 45truept{\vbox to 0pt{\vglue -13truept{\logomed A\kern -.37em{\logobig 
T}\kern -.38em G}\vss}\hss}
\break
{\small Volume \thevolumenumber\ (\thevolumeyear)
\startpage--\finishpage\nl
Published: \publishdate}

\vglue .25truein

{\parskip=0pt\leftskip 0pt plus
1fil\def\\{\par\smallskip}{\Large\bf\thetitle}\par\medskip} \vglue
0.05truein

{\parskip=0pt\leftskip 0pt plus 1fil\def\\{\par}{\sc\theauthors}
\par\medskip}%
 
\vglue 0.03truein 

{\small\leftskip 25truept\rightskip 25truept{\bf Abstract}\stdspace\theabstract

{\bf AMS Classification}\stdspace\theprimaryclass
\ifx\thesecondaryclass\relax\else; \thesecondaryclass\fi\par
{\bf Keywords}\stdspace \thekeywords\par}\vglue 7truept

}   

\ifplaintex
\font\phead=cmsl9 scaled 950
\font\pnum=cmbx10 scaled 913
\font\pfoot=cmsl9 scaled 950
\headline{\vbox to 0pt{\vskip -4.5mm\line{\small\phead\ifnum
\count0=\startpage ISSN 1472-2739 (on-line) 1472-2747 (printed)
\hfill {\pnum\folio}\else\ifodd\count0\def\\{ }%
\ifx\theshorttitle\relax\thetitle\else\theshorttitle\fi\hfill{\pnum\folio}
\else\def\\{ and }{\pnum\folio}\hfill\ifx\theshortauthors\relax\theauthors
\else\theshortauthors\fi\fi\fi}\vss}}
\footline{\vbox to 0pt{\vglue 0mm\line{\small\pfoot\ifnum\count0=\startpage
\copyright\ \gtp\hfill\else
\agt, Volume \thevolumenumber\ (\thevolumeyear)\hfill\fi}\vss}}
\else
\font=cmsl9 scaled 1050
\font\lnum=cmbx10 
\font=cmsl9 scaled 1050
\makeatletter
\def\@oddhead{{\small\lhead\ifnum\count0=\startpage ISSN 1472-2739 
(on-line) 1472-2747 (printed)\hfill {\lnum\number\count0}\else\ifodd\count0
\def\\{ }\ifx\theshorttitle\relax \thetitle \else\theshorttitle\fi\hfill
{\lnum\number\count0}\else\def\\{ and }{\lnum\number\count0}
\hfill\ifx\theshortauthors\relax 
\theauthors\else\theshortauthors\fi\fi\fi}}\def\@evenhead{\@oddhead}
\def\@oddfoot{\small\lfoot\ifnum\count0=\startpage\copyright\ \gtp\hfill\else
\agt, Volume \thevolumenumber\ (\thevolumeyear)\hfill\fi}
\def\@evenfoot{\@oddfoot}
\makeatother
\fi
\let\maketitlepage\makeagttitle
\let\makeshorttitle\maketitlepage
\let\maketitle\maketitlepage

\newwrite\gtoutfile
\long\gdef\makeheadfile{  
{\def\\{, }\def\s{ }
\immediate\openout\gtoutfile head.xxx
\immediate\write\gtoutfile{To: math@arxiv.org}
\immediate\write\gtoutfile{Subject: put OR rep NNNNN:ppppp}
\immediate\write\gtoutfile{--text follows this line--}
\immediate\write\gtoutfile{Proxy-for: \ifx\theasciiauthors\relax
\theauthors\else\theasciiauthors\fi\s<\ifx\theasciiemail\relax\theemail\else\theasciiemail\fi>}
\immediate\write\gtoutfile{\noexpand\\}
\immediate\write\gtoutfile{Authors: \ifx\theasciiauthors\relax
\theauthors\else\theasciiauthors\fi}
{\def\\{ }\immediate\write\gtoutfile{Title: \ifx\theasciititle\relax
\thetitle\else\theasciititle\fi}}
\immediate\write\gtoutfile{Subj-class: GT or SG, GR etc}
\immediate\write\gtoutfile{MSC-class: \theprimaryclass\ifx\thesecondaryclass\relax\else, \thesecondaryclass\fi}
\immediate\write\gtoutfile{Journal-ref: Algebr. Geom. Topol. \thevolumenumber\s
(\thevolumeyear) \startpage-\finishpage}
\immediate\write\gtoutfile{Comments: Published by Algebraic and
Geometric Topology at}
\immediate\write\gtoutfile{\s\s\s  http://www.maths.warwick.ac.uk/agt/AGTVol\thevolumenumber/agt-\thevolumenumber-\thepapernumber.abs.html}
\immediate\write\gtoutfile{\noexpand\\}
\immediate\write\gtoutfile{}
\ifx\theasciiabstract\relax
\immediate\write\gtoutfile{\theabstract}\else
\immediate\write\gtoutfile{\theasciiabstract}\fi
\immediate\write\gtoutfile{}
\immediate\write\gtoutfile{\noexpand\\}
\immediate\write\gtoutfile{}
\immediate\closeout\gtoutfile}}  

\def\maketitlepage{\makeagttitle\makeheadfile}
\let\makeshorttitle\maketitlepage
\let\maketitle\maketitlepage

\def\ifplaintex{\expandafter\ifx\csname documentclass\endcsname\relax}

\def\gtp{{\mathsurround=0pt\it $\cal G\mskip-2mu$eometry \&\ 
$\cal T\!\!$opology $\cal P\!$ublications}}  

\def\recd{{\small Received:\qua\receiveddate\ifx\reviseddate\relax
\else\qquad Revised:\qua\reviseddate\fi\par}} 

\def\lognumber#1{\def\thelognumber{#1}}
\def\volumenumber#1{\def\thevolumenumber{#1}}
\def\volumeyear#1{\def\thevolumeyear{#1}}
\def\papernumber#1{\def\thepapernumber{#1}}
\def\pagenumbers#1#2{\def\startpage{#1}\def\finishpage{#2}}
\def\published#1{\def\publishdate{#1}}

\def\received#1{\def\receiveddate{#1}}

\def\accepted#1{\def\accepteddate{#1}}

\def\asciiaddress#1{\def\theasciiaddress{#1}}
\def\asciiemail#1{\def\theasciiemail{#1}}

\long\def\asciiabstract#1{\long\def\theasciiabstract{#1}}
\def\asciikeywords#1{\def\theasciikeywords{#1}}

\let\\\par\let\thelognumber\relax\let\thevolumenumber\relax
\let\thepapernumber\relax\let\thevolumeyear\relax\let\startpage\relax
\let\finishpage\relax\let\publishdate\relax\let\receiveddate\relax
\let\reviseddate\relax\let\accepteddate\relax\let\theasciititle\relax
\let\theasciiauthors\relax\let\theasciiaddress\relax
\let\theasciiabstract\relax\let\theasciikeywords\relax

\let\theasciiemail\relax

\ifplaintex
\font\logobig=cmssbx10 scaled 3836
\font\logomed=cmssbx10 scaled 2557
\else
\font\logobig=cmssbx10 scaled 4200
\font\logomed=cmssbx10 scaled 2800
\fi

\long\def\makeagttitle{   
\count0=\startpage
\agt\hfill      
\hbox to 45truept{\vbox to 0pt{\vglue -13truept{\logomed A\kern -.37em{\logobig 
T}\kern -.38em G}\vss}\hss}
\break
{\small Volume \thevolumenumber\ (\thevolumeyear)
\startpage--\finishpage\nl
Published: \publishdate}

\vglue .25truein

{\parskip=0pt\leftskip 0pt plus
1fil\def\\{\par\smallskip}{\Large\bf\thetitle}\par\medskip} \vglue
0.05truein

{\parskip=0pt\leftskip 0pt plus 1fil\def\\{\par}{\sc\theauthors}
\par\medskip}%
 
\vglue 0.03truein 

{\small\leftskip 25truept\rightskip 25truept{\bf Abstract}\stdspace\theabstract

{\bf AMS Classification}\stdspace\theprimaryclass
\ifx\thesecondaryclass\relax\else; \thesecondaryclass\fi\par
{\bf Keywords}\stdspace \thekeywords\par}\vglue 7truept

}   

\ifplaintex
\font\phead=cmsl9 scaled 950
\font\pnum=cmbx10 scaled 913
\font\pfoot=cmsl9 scaled 950
\headline{\vbox to 0pt{\vskip -4.5mm\line{\small\phead\ifnum
\count0=\startpage ISSN 1472-2739 (on-line) 1472-2747 (printed)
\hfill {\pnum\folio}\else\ifodd\count0\def\\{ }%
\ifx\theshorttitle\relax\thetitle\else\theshorttitle\fi\hfill{\pnum\folio}
\else\def\\{ and }{\pnum\folio}\hfill\ifx\theshortauthors\relax\theauthors
\else\theshortauthors\fi\fi\fi}\vss}}
\footline{\vbox to 0pt{\vglue 0mm\line{\small\pfoot\ifnum\count0=\startpage
\copyright\ \gtp\hfill\else
\agt, Volume \thevolumenumber\ (\thevolumeyear)\hfill\fi}\vss}}
\else
\font=cmsl9 scaled 1050
\font\lnum=cmbx10 
\font=cmsl9 scaled 1050
\makeatletter
\def\@oddhead{{\small\lhead\ifnum\count0=\startpage ISSN 1472-2739 
(on-line) 1472-2747 (printed)\hfill {\lnum\number\count0}\else\ifodd\count0
\def\\{ }\ifx\theshorttitle\relax \thetitle \else\theshorttitle\fi\hfill
{\lnum\number\count0}\else\def\\{ and }{\lnum\number\count0}
\hfill\ifx\theshortauthors\relax 
\theauthors\else\theshortauthors\fi\fi\fi}}\def\@evenhead{\@oddhead}
\def\@oddfoot{\small\lfoot\ifnum\count0=\startpage\copyright\ \gtp\hfill\else
\agt, Volume \thevolumenumber\ (\thevolumeyear)\hfill\fi}
\def\@evenfoot{\@oddfoot}
\makeatother
\fi
\let\maketitlepage\makeagttitle
\let\makeshorttitle\maketitlepage
\let\maketitle\maketitlepage

\newwrite\gtoutfile
\long\gdef\makeheadfile{  
{\def\\{, }\def\s{ }
\immediate\openout\gtoutfile head.xxx
\immediate\write\gtoutfile{To: math@arxiv.org}
\immediate\write\gtoutfile{Subject: put OR rep NNNNN:ppppp}
\immediate\write\gtoutfile{--text follows this line--}
\immediate\write\gtoutfile{Proxy-for: \ifx\theasciiauthors\relax
\theauthors\else\theasciiauthors\fi\s<\ifx\theasciiemail\relax\theemail\else\theasciiemail\fi>}
\immediate\write\gtoutfile{\noexpand\\}
\immediate\write\gtoutfile{Authors: \ifx\theasciiauthors\relax
\theauthors\else\theasciiauthors\fi}
{\def\\{ }\immediate\write\gtoutfile{Title: \ifx\theasciititle\relax
\thetitle\else\theasciititle\fi}}
\immediate\write\gtoutfile{Subj-class: GT or SG, GR etc}
\immediate\write\gtoutfile{MSC-class: \theprimaryclass\ifx\thesecondaryclass\relax\else, \thesecondaryclass\fi}
\immediate\write\gtoutfile{Journal-ref: Algebr. Geom. Topol. \thevolumenumber\s
(\thevolumeyear) \startpage-\finishpage}
\immediate\write\gtoutfile{Comments: Published by Algebraic and
Geometric Topology at}
\immediate\write\gtoutfile{\s\s\s  http://www.maths.warwick.ac.uk/agt/AGTVol\thevolumenumber/agt-\thevolumenumber-\thepapernumber.abs.html}
\immediate\write\gtoutfile{\noexpand\\}
\immediate\write\gtoutfile{}
\ifx\theasciiabstract\relax
\immediate\write\gtoutfile{\theabstract}\else
\immediate\write\gtoutfile{\theasciiabstract}\fi
\immediate\write\gtoutfile{}
\immediate\write\gtoutfile{\noexpand\\}
\immediate\write\gtoutfile{}
\immediate\closeout\gtoutfile}}  

\def\maketitlepage{\makeagttitle\makeheadfile}
\let\makeshorttitle\maketitlepage
\let\maketitle\maketitlepage

\def\ifplaintex{\expandafter\ifx\csname documentclass\endcsname\relax}

\def\gtp{{\mathsurround=0pt\it $\cal G\mskip-2mu$eometry \&\ 
$\cal T\!\!$opology $\cal P\!$ublications}}  

\def\recd{{\small Received:\qua\receiveddate\ifx\reviseddate\relax
\else\qquad Revised:\qua\reviseddate\fi\par}} 

\def\lognumber#1{\def\thelognumber{#1}}
\def\volumenumber#1{\def\thevolumenumber{#1}}
\def\volumeyear#1{\def\thevolumeyear{#1}}
\def\papernumber#1{\def\thepapernumber{#1}}
\def\pagenumbers#1#2{\def\startpage{#1}\def\finishpage{#2}}
\def\published#1{\def\publishdate{#1}}

\def\received#1{\def\receiveddate{#1}}

\def\accepted#1{\def\accepteddate{#1}}

\def\asciiaddress#1{\def\theasciiaddress{#1}}
\def\asciiemail#1{\def\theasciiemail{#1}}

\long\def\asciiabstract#1{\long\def\theasciiabstract{#1}}
\def\asciikeywords#1{\def\theasciikeywords{#1}}

\let\\\par\let\thelognumber\relax\let\thevolumenumber\relax
\let\thepapernumber\relax\let\thevolumeyear\relax\let\startpage\relax
\let\finishpage\relax\let\publishdate\relax\let\receiveddate\relax
\let\reviseddate\relax\let\accepteddate\relax\let\theasciititle\relax
\let\theasciiauthors\relax\let\theasciiaddress\relax
\let\theasciiabstract\relax\let\theasciikeywords\relax

\let\theasciiemail\relax

\ifplaintex
\font\logobig=cmssbx10 scaled 3836
\font\logomed=cmssbx10 scaled 2557
\else
\font\logobig=cmssbx10 scaled 4200
\font\logomed=cmssbx10 scaled 2800
\fi

\long\def\makeagttitle{   
\count0=\startpage
\agt\hfill      
\hbox to 45truept{\vbox to 0pt{\vglue -13truept{\logomed A\kern -.37em{\logobig 
T}\kern -.38em G}\vss}\hss}
\break
{\small Volume \thevolumenumber\ (\thevolumeyear)
\startpage--\finishpage\nl
Published: \publishdate}

\vglue .25truein

{\parskip=0pt\leftskip 0pt plus
1fil\def\\{\par\smallskip}{\Large\bf\thetitle}\par\medskip} \vglue
0.05truein

{\parskip=0pt\leftskip 0pt plus 1fil\def\\{\par}{\sc\theauthors}
\par\medskip}%
 
\vglue 0.03truein 

{\small\leftskip 25truept\rightskip 25truept{\bf Abstract}\stdspace\theabstract

{\bf AMS Classification}\stdspace\theprimaryclass
\ifx\thesecondaryclass\relax\else; \thesecondaryclass\fi\par
{\bf Keywords}\stdspace \thekeywords\par}\vglue 7truept

}   

\ifplaintex
\font\phead=cmsl9 scaled 950
\font\pnum=cmbx10 scaled 913
\font\pfoot=cmsl9 scaled 950
\headline{\vbox to 0pt{\vskip -4.5mm\line{\small\phead\ifnum
\count0=\startpage ISSN 1472-2739 (on-line) 1472-2747 (printed)
\hfill {\pnum\folio}\else\ifodd\count0\def\\{ }%
\ifx\theshorttitle\relax\thetitle\else\theshorttitle\fi\hfill{\pnum\folio}
\else\def\\{ and }{\pnum\folio}\hfill\ifx\theshortauthors\relax\theauthors
\else\theshortauthors\fi\fi\fi}\vss}}
\footline{\vbox to 0pt{\vglue 0mm\line{\small\pfoot\ifnum\count0=\startpage
\copyright\ \gtp\hfill\else
\agt, Volume \thevolumenumber\ (\thevolumeyear)\hfill\fi}\vss}}
\else
\font=cmsl9 scaled 1050
\font\lnum=cmbx10 
\font=cmsl9 scaled 1050
\makeatletter
\def\@oddhead{{\small\lhead\ifnum\count0=\startpage ISSN 1472-2739 
(on-line) 1472-2747 (printed)\hfill {\lnum\number\count0}\else\ifodd\count0
\def\\{ }\ifx\theshorttitle\relax \thetitle \else\theshorttitle\fi\hfill
{\lnum\number\count0}\else\def\\{ and }{\lnum\number\count0}
\hfill\ifx\theshortauthors\relax 
\theauthors\else\theshortauthors\fi\fi\fi}}\def\@evenhead{\@oddhead}
\def\@oddfoot{\small\lfoot\ifnum\count0=\startpage\copyright\ \gtp\hfill\else
\agt, Volume \thevolumenumber\ (\thevolumeyear)\hfill\fi}
\def\@evenfoot{\@oddfoot}
\makeatother
\fi
\let\maketitlepage\makeagttitle
\let\makeshorttitle\maketitlepage
\let\maketitle\maketitlepage

\newwrite\gtoutfile
\long\gdef\makeheadfile{  
{\def\\{, }\def\s{ }
\immediate\openout\gtoutfile head.xxx
\immediate\write\gtoutfile{To: math@arxiv.org}
\immediate\write\gtoutfile{Subject: put OR rep NNNNN:ppppp}
\immediate\write\gtoutfile{--text follows this line--}
\immediate\write\gtoutfile{Proxy-for: \ifx\theasciiauthors\relax
\theauthors\else\theasciiauthors\fi\s<\ifx\theasciiemail\relax\theemail\else\theasciiemail\fi>}
\immediate\write\gtoutfile{\noexpand\\}
\immediate\write\gtoutfile{Authors: \ifx\theasciiauthors\relax
\theauthors\else\theasciiauthors\fi}
{\def\\{ }\immediate\write\gtoutfile{Title: \ifx\theasciititle\relax
\thetitle\else\theasciititle\fi}}
\immediate\write\gtoutfile{Subj-class: GT or SG, GR etc}
\immediate\write\gtoutfile{MSC-class: \theprimaryclass\ifx\thesecondaryclass\relax\else, \thesecondaryclass\fi}
\immediate\write\gtoutfile{Journal-ref: Algebr. Geom. Topol. \thevolumenumber\s
(\thevolumeyear) \startpage-\finishpage}
\immediate\write\gtoutfile{Comments: Published by Algebraic and
Geometric Topology at}
\immediate\write\gtoutfile{\s\s\s  http://www.maths.warwick.ac.uk/agt/AGTVol\thevolumenumber/agt-\thevolumenumber-\thepapernumber.abs.html}
\immediate\write\gtoutfile{\noexpand\\}
\immediate\write\gtoutfile{}
\ifx\theasciiabstract\relax
\immediate\write\gtoutfile{\theabstract}\else
\immediate\write\gtoutfile{\theasciiabstract}\fi
\immediate\write\gtoutfile{}
\immediate\write\gtoutfile{\noexpand\\}
\immediate\write\gtoutfile{}
\immediate\closeout\gtoutfile}}  

\def\maketitlepage{\makeagttitle\makeheadfile}
\let\makeshorttitle\maketitlepage
\let\maketitle\maketitlepage

\def\ifplaintex{\expandafter\ifx\csname documentclass\endcsname\relax}

\def\gtp{{\mathsurround=0pt\it $\cal G\mskip-2mu$eometry \&\ 
$\cal T\!\!$opology $\cal P\!$ublications}}  

\def\recd{{\small Received:\qua\receiveddate\ifx\reviseddate\relax
\else\qquad Revised:\qua\reviseddate\fi\par}} 

\def\lognumber#1{\def\thelognumber{#1}}
\def\volumenumber#1{\def\thevolumenumber{#1}}
\def\volumeyear#1{\def\thevolumeyear{#1}}
\def\papernumber#1{\def\thepapernumber{#1}}
\def\pagenumbers#1#2{\def\startpage{#1}\def\finishpage{#2}}
\def\published#1{\def\publishdate{#1}}

\def\received#1{\def\receiveddate{#1}}

\def\accepted#1{\def\accepteddate{#1}}

\def\asciiaddress#1{\def\theasciiaddress{#1}}
\def\asciiemail#1{\def\theasciiemail{#1}}

\long\def\asciiabstract#1{\long\def\theasciiabstract{#1}}
\def\asciikeywords#1{\def\theasciikeywords{#1}}

\let\\\par\let\thelognumber\relax\let\thevolumenumber\relax
\let\thepapernumber\relax\let\thevolumeyear\relax\let\startpage\relax
\let\finishpage\relax\let\publishdate\relax\let\receiveddate\relax
\let\reviseddate\relax\let\accepteddate\relax\let\theasciititle\relax
\let\theasciiauthors\relax\let\theasciiaddress\relax
\let\theasciiabstract\relax\let\theasciikeywords\relax

\let\theasciiemail\relax

\ifplaintex
\font\logobig=cmssbx10 scaled 3836
\font\logomed=cmssbx10 scaled 2557
\else
\font\logobig=cmssbx10 scaled 4200
\font\logomed=cmssbx10 scaled 2800
\fi

\long\def\makeagttitle{   
\count0=\startpage
\agt\hfill      
\hbox to 45truept{\vbox to 0pt{\vglue -13truept{\logomed A\kern -.37em{\logobig 
T}\kern -.38em G}\vss}\hss}
\break
{\small Volume \thevolumenumber\ (\thevolumeyear)
\startpage--\finishpage\nl
Published: \publishdate}

\vglue .25truein

{\parskip=0pt\leftskip 0pt plus
1fil\def\\{\par\smallskip}{\Large\bf\thetitle}\par\medskip} \vglue
0.05truein

{\parskip=0pt\leftskip 0pt plus 1fil\def\\{\par}{\sc\theauthors}
\par\medskip}%
 
\vglue 0.03truein 

{\small\leftskip 25truept\rightskip 25truept{\bf Abstract}\stdspace\theabstract

{\bf AMS Classification}\stdspace\theprimaryclass
\ifx\thesecondaryclass\relax\else; \thesecondaryclass\fi\par
{\bf Keywords}\stdspace \thekeywords\par}\vglue 7truept

}   

\ifplaintex
\font\phead=cmsl9 scaled 950
\font\pnum=cmbx10 scaled 913
\font\pfoot=cmsl9 scaled 950
\headline{\vbox to 0pt{\vskip -4.5mm\line{\small\phead\ifnum
\count0=\startpage ISSN 1472-2739 (on-line) 1472-2747 (printed)
\hfill {\pnum\folio}\else\ifodd\count0\def\\{ }%
\ifx\theshorttitle\relax\thetitle\else\theshorttitle\fi\hfill{\pnum\folio}
\else\def\\{ and }{\pnum\folio}\hfill\ifx\theshortauthors\relax\theauthors
\else\theshortauthors\fi\fi\fi}\vss}}
\footline{\vbox to 0pt{\vglue 0mm\line{\small\pfoot\ifnum\count0=\startpage
\copyright\ \gtp\hfill\else
\agt, Volume \thevolumenumber\ (\thevolumeyear)\hfill\fi}\vss}}
\else
\font=cmsl9 scaled 1050
\font\lnum=cmbx10 
\font=cmsl9 scaled 1050
\makeatletter
\def\@oddhead{{\small\lhead\ifnum\count0=\startpage ISSN 1472-2739 
(on-line) 1472-2747 (printed)\hfill {\lnum\number\count0}\else\ifodd\count0
\def\\{ }\ifx\theshorttitle\relax \thetitle \else\theshorttitle\fi\hfill
{\lnum\number\count0}\else\def\\{ and }{\lnum\number\count0}
\hfill\ifx\theshortauthors\relax 
\theauthors\else\theshortauthors\fi\fi\fi}}\def\@evenhead{\@oddhead}
\def\@oddfoot{\small\lfoot\ifnum\count0=\startpage\copyright\ \gtp\hfill\else
\agt, Volume \thevolumenumber\ (\thevolumeyear)\hfill\fi}
\def\@evenfoot{\@oddfoot}
\makeatother
\fi
\let\maketitlepage\makeagttitle
\let\makeshorttitle\maketitlepage
\let\maketitle\maketitlepage

\newwrite\gtoutfile
\long\gdef\makeheadfile{  
{\def\\{, }\def\s{ }
\immediate\openout\gtoutfile head.xxx
\immediate\write\gtoutfile{To: math@arxiv.org}
\immediate\write\gtoutfile{Subject: put OR rep NNNNN:ppppp}
\immediate\write\gtoutfile{--text follows this line--}
\immediate\write\gtoutfile{Proxy-for: \ifx\theasciiauthors\relax
\theauthors\else\theasciiauthors\fi\s<\ifx\theasciiemail\relax\theemail\else\theasciiemail\fi>}
\immediate\write\gtoutfile{\noexpand\\}
\immediate\write\gtoutfile{Authors: \ifx\theasciiauthors\relax
\theauthors\else\theasciiauthors\fi}
{\def\\{ }\immediate\write\gtoutfile{Title: \ifx\theasciititle\relax
\thetitle\else\theasciititle\fi}}
\immediate\write\gtoutfile{Subj-class: GT or SG, GR etc}
\immediate\write\gtoutfile{MSC-class: \theprimaryclass\ifx\thesecondaryclass\relax\else, \thesecondaryclass\fi}
\immediate\write\gtoutfile{Journal-ref: Algebr. Geom. Topol. \thevolumenumber\s
(\thevolumeyear) \startpage-\finishpage}
\immediate\write\gtoutfile{Comments: Published by Algebraic and
Geometric Topology at}
\immediate\write\gtoutfile{\s\s\s  http://www.maths.warwick.ac.uk/agt/AGTVol\thevolumenumber/agt-\thevolumenumber-\thepapernumber.abs.html}
\immediate\write\gtoutfile{\noexpand\\}
\immediate\write\gtoutfile{}
\ifx\theasciiabstract\relax
\immediate\write\gtoutfile{\theabstract}\else
\immediate\write\gtoutfile{\theasciiabstract}\fi
\immediate\write\gtoutfile{}
\immediate\write\gtoutfile{\noexpand\\}
\immediate\write\gtoutfile{}
\immediate\closeout\gtoutfile}}  

\def\maketitlepage{\makeagttitle\makeheadfile}
\let\makeshorttitle\maketitlepage
\let\maketitle\maketitlepage

\lognumber{14}
\volumenumber{2}
\volumeyear{2002}
\papernumber{14}
\published{24 April 2002}
\pagenumbers{297}{309}
\received{7 January 2002}
\accepted{6 March 2002}

\usepackage{amssymb, amsmath}
\usepackage{epsf}

\newtheorem{pro}{Proposition}[section]
\newtheorem{thm}[pro]{Theorem}
\newtheorem{lem}[pro]{Lemma}

\newtheorem*{qthm}{Theorem}
\theoremstyle{definition}
\newtheorem{dfn}[pro]{Definition}

\theoremstyle{remark}
\newtheorem*{rmk}{Remark}

\newcommand{\s}{\Sigma}
\newcommand{\ie}{ie,}
\newcommand{\eg}{eg,}

\newcommand{\Cf}{Cf}

\newcommand{\del}{\partial}

\begin{document}

\title{Thin position for a connected sum of small knots}

\authors{Yo'av Rieck\\Eric Sedgwick}
\address{Department of Mathematics, Nara Women's University\\
       Kitauoya Nishimachi, Nara 630-8506, Japan, and\\
       Department of Mathematics, University of Arkansas\\
       Fayetteville, AR 72701, USA}
\email{yoav@mail.uark.edu}
\secondaddress{DePaul University,
Department of Computer Science\\243 S. Wabash Ave. - Suite 401,
Chicago, IL 60604, USA} 

\secondemail{esedgwick@cs.depaul.edu}

\asciiaddress{Department of Mathematics, Nara Women's University\\
Kitauoya Nishimachi, Nara 630-8506, Japan, and\\
Department of Mathematics, University of Arkansas\\
Fayetteville, AR 72701, USA\\and\\DePaul University,
Department of Computer Science\\243 S. Wabash Ave. - Suite 401,
Chicago, IL 60604, USA}
\asciiemail{yoav@mail.uark.edu, esedgwick@cs.depaul.edu}

\begin{abstract}
We show that every thin position for a connected sum of small
knots is obtained in an obvious way:  place each summand in thin
position so that no two summands intersect the same level surface,
then connect the lowest minimum of each summand to the highest
maximum of the adjacent summand below.  See Figure
\ref{fig:thin-position}.
\end{abstract}
\asciiabstract{
We show that every thin position for a connected sum of small
knots is obtained in an obvious way:  place each summand in thin
position so that no two summands intersect the same level surface,
then connect the lowest minimum of each summand to the highest
maximum of the adjacent summand below.}


\asciikeywords{3-manifold, connected sum of knots, thin position}

\maketitle

\section{Introduction}

In \cite{gabai} Gabai introduced {\it thin position}, a tool since
used for many results on knots and graphs in $S^3$: Gabai's proof
of Property $R$ \cite{gabai}, Gordon and Luecke's \cite{g-l} proof
of the Knot Complement Theorem, and Scharlemann and Thompson's
\cite{s-t} proof of Waldhausen's Theorem that irreducible Heegaard
splittings of $S^3$ are unique, among others. Although thin
position has been used in more general settings (\eg\ by Rieck
\cite{rieck:topology} and Rieck and Sedgwick \cite{rs.daisy} to
study the behavior of Heegaard Surfaces under Dehn surgery) it has
been most fruitful when applied to the study of knots in $S^3$. We
examine thin position in $S^3$ for a connected sum of knots.

\begin{figure}[ht!]
{\epsfxsize = 3 in \centerline{\epsfbox{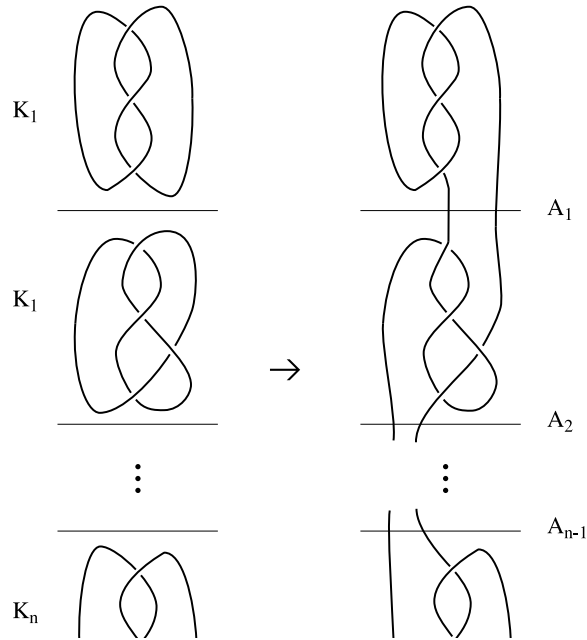}} }
\label{fig:thin-position} \caption{Thin position for a connected
sum of small knots is a stack of the summands.}
\end{figure}

Given  a connected sum of knots, say $K=K_1\#..\#K_n$, there are
obvious candidates for thin position of $K$.  Choose an ordering
of the summands $K_1, ..., K_n$.   and put each in a thin position
so that successive pairs are separated by level $2$-spheres.  Form
the connected sum $K$ by connecting the lowest minimum of $K_i$ to
the highest maximum of $K_{i+1}$.  Each level sphere is punctured
twice and becomes a decomposing annulus which appears in this
presentation as a {\it thin level} surface, a surface with a
minimum of the knot immediately above and a maximum of the knot
immediately below.   See Figure \ref{fig:thin-position}. We call
this position for $K$ a \it stack \rm on the summands $K_1, ...,
K_n$. The width of any stack on these summands is $\Sigma_{i=1}^n
w(K_i)-2(n-1)$, which gives an upper bound for the width of $K$.

Is this thin position for $K$?  In fact, it is not immediately
obvious that thin position for a connected sum has any thin levels
whatsoever, thin position could be bridge position.  While
Thompson has shown that small knots do not have thin levels
\cite{thompson:ThinBridge}, the converse is simply false. In
\cite{heath-kobayashi} Heath and Kobayashi give an example of a
knot (originally considered by Morimoto \cite{morimoto}) for which
thin position is bridge position even though its exterior contains
an essential four times punctured sphere.  Our first theorem gives
a converse to Thompson's theorem in the case of a connected sum:

\begin{qthm}[\ref{thm:TPneqB}]
Let $K$ be a connected sum of non-trivial knots, $ K\! =\! K_1 \#
K_2$. Then thin position for $K$ is not bridge position for $K$.
\end{qthm}

We then analyze thin levels.  For $mp$--small knots (inclusive of
small knots, see Section \ref{sec:prelims} for a definition), we
demonstrate that each thin level is a decomposing annulus. Indeed,
stacks are thin:\eject

\begin{qthm}[\ref{thm:connect-small-knots}]
Let $K=\#_{i=1}^n K_i$ be a connected sum of $mp$--small knots. If
$K$ is in thin position, then there is an ordering of the summands
$K_{i_1}, K_{i_2},...,K_{i_n}$ and a collection of leveled
decomposing annuli $A_{i_1}, A_{i_2}, ..., A_{i_{n-1}}$ so that
the thin levels of the presentation are precisely the annuli
$\{A_{i_j}\}$ occurring in order, where the annulus $A_{i_j}$
separates the connected sum $K_{i_1}\#K_{i_2}\#...\#K_{i_j}$ from
the connected sum $K_{i_{j+1}}\#...\#K_{i_n}$. {\rm (See Figure
\ref{fig:thin-position}.)}
\end{qthm}

For $mp$--small knots this yields the equality $w(K) = \Sigma_{i=1}^n
w(K_i) -2(n-1)$.  Finally, since $mp$--small knots in thin position
are also in bridge position \cite{thompson:ThinBridge}, each
component of the stack is in bridge position.  After connecting
these knots, the number of maxima is $\Sigma_{i=1}^n b(K_i) -
(n-1)$. Thus, while thin position is not bridge position for these
knots, it realizes the minimal bridge number as given by Schubert
\cite {schubert} (see Schultens \cite{schultens:bridge} for a
modern proof).

\rk{Acknowledgements} The authors would like to thank Marc
Lackenby and Tsuyoshi Kobayashi for helpful conversations, and
RIMS of Kyoto University and Nara Women's University for their
kind hospitality.  The first named author was supported by JSPS
grant number P00024.

\section{Preliminaries}
\label{sec:prelims}

Most of the definitions we use are standard; however we need:

\begin{dfn}
\label{dfn:mps} A knot $K \subset S^3$ is called {\it meridionally
planar small} if there is no meridional essential planar surface
in its complement.  We use the notation {\it $mp$--small}.
\end{dfn}

The set of knots under consideration is substantial and inclusive
of small knots: by CGLS \cite{cgls} small knots in $S^3$ are
meridionally small, and by definition meridionally small knots are
$mp$--small.

The width of a presentation of a link $L$ (\ie\ the width of an
embedding of a compact 1-manifold into $S^3$) is defined as
follows: let $h:S^3 \to \mathbb [-\infty,\infty]$ be a height
function that sends one point to $-\infty$, another to $\infty$
and has all other level sets diffeomorphic to $S^2$ (\ie\ a Morse
function with one minimum, one maximum and no other critical
points). Suppose $h|_L$ is Morse (else the width is not defined).
Pick a regular value between every two consecutive critical levels
of $h|_L$ and count the number of times $L$ intersects that level.
The sum of these numbers is the {\it width of the presentation}.
{\it Thin position} for $L$ is any embedding ambient isotopic to
$L$ that minimizes the width in that isotopy class.  The {\it
width of the link $L$}, denoted $w(L)$, is the width of a thin
position for $L$.

A {\it thin level} is a level set of a regular value so that the
first critical point of $h|_L$ above is a minimum and the first
critical point below is a maximum. An embedding has no thin levels
if and only if each maximum of $L$ occurs above each minimum of
$L$.   Among all embeddings of $L$ without thin levels, one with a
minimal number of maxima (equivalently minima) is called a {\it
bridge  position for $L$}.  The {\it bridge number of $L$},
denoted $b(L)$, is the number of maxima (equivalently minima) in a
bridge position for $L$.  By pulling maxima up and minima down, it is
easy to see that this is the minimal number of maxima in the isotopy
class of the knot.

\begin{figure}[ht!]
    {\epsfxsize = 3 in \centerline{\epsfbox{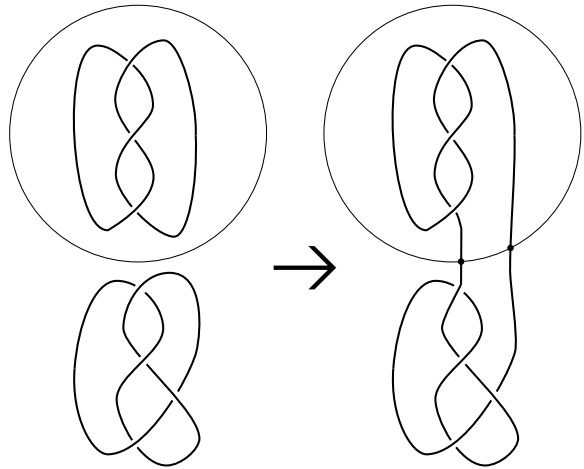}} }
    \label{fig:connect-sum}
    \caption{Connected sum of knots.}
\end{figure}

Given two knots their connected sum is described schematically in
Figure \ref{fig:connect-sum}.  For more detail, see
\cite{burde-zie}.

We recall that every knot has a unique representation as a
connected sum of prime knots, \ie\ any knot $K$ can be written
uniquely up to reordering as $\#_{i=1}^n K_i$ for some $K_i$, and
$K_i$ contains no meridional essential embedded annuli in its
exterior. Since any meridional essential embedded annulus
decomposes the knot as a connected sum (the factors are not
necessarily prime) we call such annulus a {\it decomposing
annulus}.  Note that the uniqueness referred to above is
uniqueness of factors, not of the decomposing annuli.  A {\it
split link} on knots $K_1$ and $K_2$ is an embedding of the two
knots that can be separated by an embedded $S^2$.

\section{\!Essential meridional planar surfaces in a
connected sum of $mp$--small knots}


\begin{lem}
\label{thm:EssMPsurfaces} Let $K = \#_{i=1}^n K_i$ be a connected
sum of $mp$--small knots $K_i,\; i=1, \dots ,n$. Then any essential
meridional planar surface in the exterior of $K$ is a decomposing
annulus.
\end{lem}

\begin{proof} By way of contradiction, let $n \geq 1$ be the smallest $n$
for which there is a connected sum with $n$ components that is a
counterexample to lemma.  Let $P$ be a non-annular planar essential
surface in the exterior of $K$. By assumption $n>1$, so we can choose an
annulus $A$ that decomposes $K$ as a connected sum $K=K_1 \# K_2$. The
annulus $A$ is properly embedded in $X = S^3 \setminus N(K)$ and divides
$X$ into manifolds $X_1$ and $X_2$, the exterior of $K_1$ and $K_2$
respectively.

Since $A$ and $P$ are essential surfaces in $X$, $P$ may be
isotoped to intersect $A$ essentially.  Since both surfaces have
meridional boundary, the intersection will be a (perhaps empty)
collection of simple closed curves, each curve essential in both
surfaces. Among all such positions for $P$, choose one that minimizes
the number of curves in the intersection $A \cap P$. Let $P_i = P
\cap X_i,\; i=1,2$.
Each component of the surface $P_i$ is
properly embedded in the knot exterior $X_i$ and is planar. Since
$P$ is essential and each curve of $P \cap A$ is essential in $P$,
each component of $P_i$ is necessarily incompressible in $X_i,\;
i=1,2$.

In knot exteriors ($X_1$ and $X_2$ in our case) boundary
compression implies compression for any
surface that is not an annulus.  If either $P_1$ or $P_2$ contains
a component that is not an annulus, then that component is also
boundary incompressible. This yields an essential surface which is
not an annulus, with its boundary on $A$, a meridional slope. Such
a component would contradict our assumption that $n$ was the least
$n$ so that a connected sum of $n$ components contained an
essential planar surface with meridional slope.

We conclude that each $P_i$ consists entirely of annuli.
This contradicts our assumption that $P$ was a planar surface
that is not an annulus.
\end{proof}

\section{Thin position is not bridge position}

\begin{thm}
\label{thm:TPneqB} Let $K$ be a connected sum of non-trivial
knots, $ K = K_1 \# K_2$. Then any thin position for $K$ is not bridge
position for $K$.
\end{thm}

\begin{proof}
Let $K$ be $K_1 \# K_2$, where $K_1$ and $K_2$ are not unknots. We
will prove the theorem by showing that $K$ has a position that has
lower width than its bridge position. Denote the bridge number of
$K$ by $b$ and that of $K_i$ by $b_i$ ($i=1,2$). By Schubert's
Theorem (\cite{schubert}, see \cite{schultens:bridge} for a modern
proof) $b = b_1 + b_2 - 1$.

For a knot in bridge position with $n$ maxima the width of the
presentation is easily seen to be
\begin{eqnarray}
w(K)&=& 2+4+...+(2n-2)+2n+(2n-2)+... +4+2 \nonumber \\ &=&
2(\Sigma_{i=1}^n 2i) - 2n \nonumber \\ & =&  4\frac{n(n+1)}{2} -
2n \nonumber \\& = & 2n^2. \nonumber
\end{eqnarray}
Applying this to $K$, all we need is to show that $2b^2$ is not
the minimal width for $K$.

Consider a presentation of the split link of $K_1$ and $K_2$ with
$K_1$ above the level 0 and $K_2$ below the level 0.  If we put
both knots in bridge position, the width of this presentation is
$2b_1^2 + 2b_2^2$. We obtain a presentation of $K=K_1 \# K_2$ by
connecting the lowest minimum of $K_1$ to the highest maximum of
$K_2$ without introducing new critical points, which lowers the
width by 2. Thus, it suffices to show that $2b^2 - (2b_1^2 +
2b_2^2 - 2) > 0$. Applying Schubert's Theorem we get:
\begin{eqnarray}
2b^2 - (2b_1^2 + 2b_2^2 - 2) & = &2(b_1 + b_2 - 1)^2 - (2b_1^2 +
2b_2^2 -2) \nonumber \\  & = &
[2b_1^2+2b_2^2+4b_1b_2-4(b_1+b_2)+2] - [2b_1^2 + 2b_2^2 - 2]
\nonumber
\\ & = & 4b_1b_2 - 4(b_1+b_2)+4. \nonumber
\\ & = & 4 [b_1(b_2-1) - (b_2-1)] \nonumber
\\ & = & 4 (b_1-1)(b_2-1). \nonumber
\end{eqnarray}
Since $K_1$ and $K_2$ are both non-trivial knots, $b_1 > 1$ and
$b_2 > 1$ so the above product is positive.
\end{proof}

\section{Connected sum of small knots}
\label{sec:connect-small-knots}

\begin{thm}
\label{thm:connect-small-knots} Let $K=\#_{i=1}^n K_i$ be a
connected sum of $mp$--small knots. If $K$ is in thin position, then
there is an ordering of the summands $K_{i_1},
K_{i_2},...,K_{i_n}$ and a collection of decomposing annuli
$A_{i_1}, A_{i_2}, ..., A_{i_{n-1}}$ so that the thin levels of
the presentation are precisely the annuli $\{A_{i_j}\}$ occurring
in order, where the annulus $A_{i_j}$ separates the connected sum
$K_{i_1}\#K_{i_2}\#...\#K_{i_j}$ from the connected sum
$K_{i_{j+1}}\#...\#K_{i_n}$. {\rm(See Figure \ref{fig:thin-position}.)}
\end{thm}

\begin{proof}
We induct on $n$.  Our goal is to show that there exists a leveled
decomposing annulus, \ie\ an essential annulus that is $h^{-1}(t)$
for some $t$.  The theorem will follow since above (and below) the
leveled annulus we see summands of $K$ in thin position (else we
can reduce the width of $K$).  Since the prime summands are
$mp$--small, Thompson's result \cite{thompson:ThinBridge} guarantees
that their thin positions are bridge positions, thus there are no
thin level surfaces other than the decomposing annuli.

Existence of a leveled annulus follows in two steps given by the
lemmata below.

\begin{lem}
\label{lem:arc-on-thin-level} Let $K$ be the connected sum of
$mp$--small knots and $T$ be any thin level in a thin presentation of
$K$. Then a neighborhood of a spanning arc of a decomposing
annulus is isotopic onto $T$ by an isotopy preserving $K$ setwise.
\end{lem}

\begin{figure}[ht!]
    {\epsfxsize = 2.5 in \centerline{\epsfbox{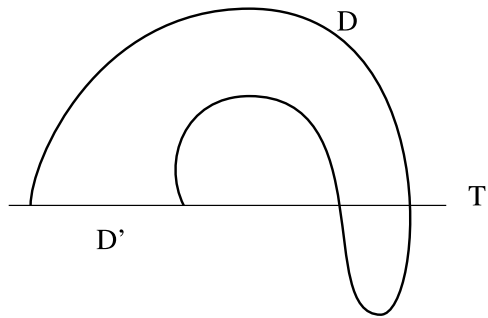}} }
    \caption{$D'$ is incompressible.}
    \label{fig:CompressT}
\end{figure}

The existence of a thin level is guaranteed by Theorem
\ref{thm:TPneqB}.

\begin{lem}
\label{lem:level-is-annulus} Let $K$ be a connected sum in thin
position, and $T$ a thin level onto which a spanning arc of a
decomposing annulus is isotopic. Then $T$ is a decomposing
annulus.
\end{lem}

\begin{rmk}
Lemma \ref{lem:level-is-annulus} does not require the assumption
that the knots are $mp$--small.  This assumption is only used when,
in the proof of Lemma \ref{lem:arc-on-thin-level}, we claim that
the essential surface we find is (some) decomposing annulus.
\end{rmk}

We may assume that there is a thin level $T$ that is inessential,
for otherwise both lemmata follow from Lemma \ref{thm:EssMPsurfaces}.

\begin{proof}[Proof of Lemma \ref{lem:arc-on-thin-level}](\Cf\
Thompson \cite{thompson:ThinBridge}.)\qua Let $T=h^{-1}(t_0)$ be a
thin level.  We assumed that $T$ is compressible.  Let $D$ be a
compressing disk for $T$, and $D'$ a punctured disk that $\del D$
bounds on $T$. Passing to an innermost compression on $D'$ we
replace $D'$ by the innermost disk and accordingly change $D$.
(Note: we may not assume that $D$ is a compressing disk for $T$ as
it may cross the level $T$ in its interior, see Figure
\ref{fig:CompressT}; we allow $D$ to compress $D'$ from above or
from below).

 We show that $D' \cup D$ is essential.  If not, it either
compresses or boundary compresses.  That $D' \cup D$ does not
compress follows from our innermost choice of $D'$ and the fact
that the boundary of any compressing disk can be isotoped to be
disjoint from $D$. The boundary of any boundary compressing disk
that lies in $D' \cup D$ can be isotoped to lie entirely in $D'$
and is hence a boundary compression for the thin level $T$.  A
boundary compression that joins a boundary component of $D'$ to
itself easily implies a compression, which we have already noted
does not occur. And, a boundary compression that connects two
distinct boundary components can never occur on a thin level. Say
that one does, and that it starts above the thin level $T$.  By
definition of thin level, the first critical point on the knot
above $T$ is a minimum. The boundary compressing disk can be used
to isotope an arc of the knot to lie below $T$ with just a single
maximum. This either pulls a maximum on the arc below the minimum
lying above $T$, or eliminates both. Either reduces the width of
the presentation. If the arc contains additional critical points
of the knot, they are eliminated in pairs, further reducing the
width.

Thus, by Theorem \ref{thm:EssMPsurfaces}, $D' \cup D$ is a
decomposing annulus.  A spanning arc for this annulus
can be isotoped to be disjoint from $D$
and therefore this arc and its neighborhood lie in $D'$, hence in
$T$.
\end{proof}

\begin{proof}[Proof of Lemma \ref{lem:level-is-annulus}]
This lemma follows from a width calculation. By\break Lem\-ma
\ref{lem:arc-on-thin-level} when we cut $K$ open along a leveled
spanning arc on $T$ we get a decomposition of $K$ into $K_+$ and
$K_-$, where $K_+$ and $K_-$ are non-trivial knots.
See Figure \ref{fig:decomposingK}.

\begin{figure}[ht!]
    {\epsfxsize = 4.25 in \centerline{\epsfbox{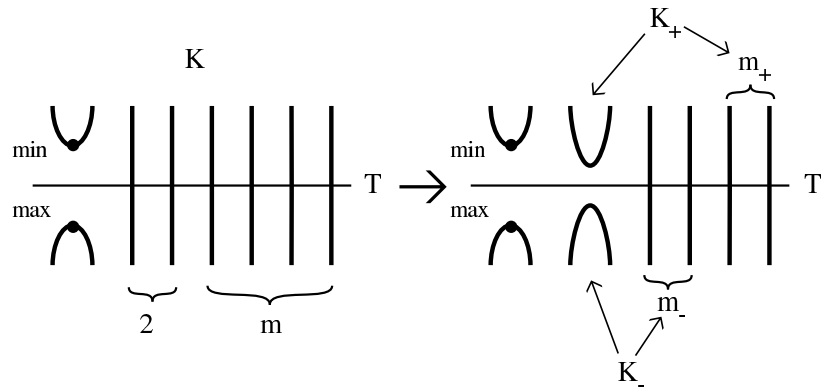}} }
    \caption{Splitting $K$ to form $K_+$ and $K_-$.}
    \label{fig:decomposingK}
\end{figure}

\begin{figure}[ht!]
    {\epsfxsize = 5 in \centerline{\epsfbox{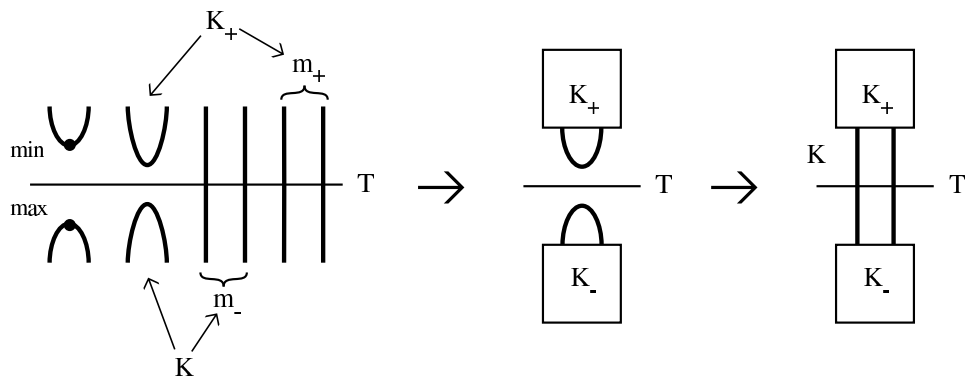}} }
    \caption{Separating $K_+$ and $K_-$ and reconnecting to obtain $K$.}
    \label{fig:Thinner}
\end{figure}

The lemma (and theorem) would follow once we show that $K_+ \cap
T$ and $K_- \cap T$ are both empty.  Denote the number of times
$K_+$ intersect $T$ by $m_+$, and the number of times $K_-$
intersect $T$ by $m_-$. Denoting $m=m_++m_-$, our goal is to show
that $m=0$.   We do this by separating $K_+$ and $K_-$ and
reconnecting to obtain a new presentation of $K$.  See Figure
\ref{fig:Thinner}.  If $m>0$ we will demonstrate that this
manipulation reduces width.

Since $T$ is a thin level for $K$ the first critical point above
it is a minimum and the first critical point below it is a
maximum. Denote these points by {\tt min} and {\tt max}, see
Figure \ref{fig:decomposingK}. After cutting $K$ open there is
another minimum above $T$, which belongs to $K_+$, and another
maximum below $T$ belonging to $K_-$.  (We do not know to which of
the knots {\tt max} and {\tt min} belong.)

We pair each maximum of $K_+$ with a minimum of $K_+$ so that for each
max\-imum-minimum pair the minimum is below the maximum. That such
correspondence exists follows easily from that fact that above each level
there are no less maxima than minima: label the maxima and label the
minima, both in descending order, and pair the $i^{th}$ maximum with the
$i^{th}$ minimum, for example. Similarly we pair up maxima and minima of
$K_-$.

The following are well known and easy facts about the calculus of
width. While isotoping a link, if two maxima change relative
heights the width does not change, and similarly for minima.
However, if we move a maximum above a minimum the width is raised
by four, while isotoping a minimum above a maximum the width is
lowered by four, see Figure \ref{fig:move-rel-extrema}.

\begin{figure}[ht!]
    {\epsfxsize = \hsize \centerline{\epsfbox{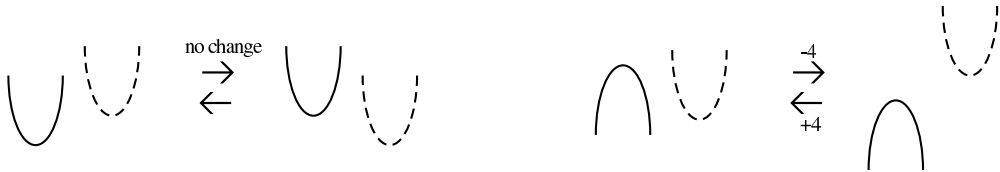}} }
    \caption{Exchanging relative extrema.}
    \label{fig:move-rel-extrema}
\end{figure}

We now begin our width calculation: our first move was cutting $K$
and obtaining the split link on $K_+$ and $K_-$.  By doing so, we
removed the level $T$ (with $m+2$ punctures) and replaced it by
three critical levels, two of width $m+2$ and one of width $m$.
Thus we raised the width by $2m+2$.  Next we isotope $K_+$ rigidly
to lie above $T$ and $K_-$ to lie below it.  See Figure
\ref{fig:Thinner}.  During this isotopy $K_+$ and $K_-$ may
intersect each other, but by the end of the process we will once
again have the split link on $K_+$ and $K_-$. The isotopy may also
temporarily increase the width as a maximum of $K_+$ is isotoped
past a minimum of $K_-$; however that contribution will be
canceled when the corresponding minimum of $K_+$ is isotoped past
the corresponding maximum of $K_-$. Finally, we will connect the
lowest minimum of $K_+$ to the highest maximum of $K_-$, obtaining
again a presentation of $K$.

For the next definition, see Figure \ref{fig:split-pair}.

\begin{dfn} \label{dfn:split-pair}
Let $X_+,\;x_+\;Y_-,\;y_-$ be four critical points so that $X_+$
(resp. $Y_-$) is a maximum of $K_+$ (resp. $K_-$) and $x_+$ (resp.
$y_-$) is a minimum of $K_+$ (resp. $K_-$). Assume further that
$X_+$ is paired with $x_+$ and $Y_-$ is paired with $y_-$.  Then
the pair $((X_+,x_+),(Y_-,y_-))$ is called a {\it split pair} if
$h(X_+)>h(y_-)$ and $h(x_+) < h(Y_-)$.
\end{dfn}

\begin{figure}[ht!]
    {\epsfxsize = 4 in \centerline{\epsfbox{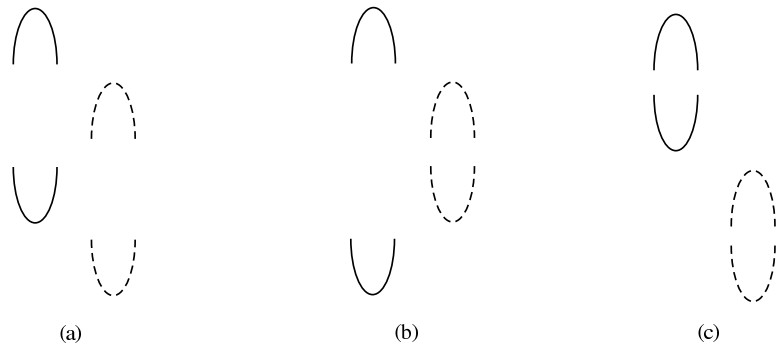} }}
    \caption{Two split pairs and a non-split pair.}
    \label{fig:split-pair}
\end{figure}

Split pairs are exactly those pairs which lower the width when we
separate $K_+$ and $K_-$: the maximum $X_+$ is already higher than
the minimum $y_-$, so the width is not raised, but the minimum
$x_+$ is below the maximum $Y_-$, so the width is lowered by four.
In Figure \ref{fig:split-pair} one of the of knots (either dashed
or solid) is $K_+$ and the other is $K_-$. Figures (a) and (b) are
both split pairs and (c) is not, independent of which knot we
choose as $K_+$ and which $K_-$.  Being a split pair is a
geometric property: both maxima have to be above both minima.
Separating a split pair (by moving one knot up and the other down)
reduces the width by four regardless of the direction we move the
knots.

To get a lower bound on the width reduction we achieve, we must
estimate the number of split pairs. There are four possibilities:

\begin{enumerate}
\item $\mbox{\tt min} \in K_+$ and $\mbox{\tt max} \in K_-$
\item $\mbox{\tt min} \in K_-$ and $\mbox{\tt max} \in K_+$
\item $\mbox{\tt min} \in K_+$ and $\mbox{\tt max} \in K_+$
\item $\mbox{\tt min} \in K_-$ and $\mbox{\tt max} \in K_-$
\end{enumerate}

$K_+$ has exactly $\frac{m_+}{2}$ maximum-minimum pairs in which
the maximum is above the level $T$ ($\frac{m_+}{2}$ is the number
of arcs $K_+$ has in the ball above $T$ and in the ball below it);
similarly $K_-$ has exactly $\frac{m_-}{2}$ pairs separated by
$T$.  Since all of these minima are below $T$ and the maxima
above, we note that the two minima adjacent to $T$ above it and
the two maxima below it are not members of these pairs.  (We make
this comment to ensure that no split pair is counted twice.)  With
this in mind we are ready to treat each of the cases above:

\begin{enumerate}
\item In this case $K_+$ has two minima directly above $T$.
Each of these minima (together with its corresponding maximum)
will be involved in a split pair with each of the pairs of $K_-$
separated by $T$, a total of $\frac{m_-}{2}$ split pairs.
Similarly, the two maxima of $K_-$ below $T$ will be involved in
$\frac{m_+}{2}$ split pairs each.  Since the minima above $T$ and
the maxima below it are not members of pairs separated by $T$ no
split pair is counted twice. The width is lowered by four per
split pair, yielding a reduction of at least $4(m_++m_-)=4m$.

\item  As in Case (1) the minimum directly above $T$, which is a point of
$K_+$, is involved in $\frac{m_-}{2}$ split pairs.  But it is also
involved in a split pair with {\tt min} and its corresponding
maximum, and is thus involved in a total of $\frac{m_-}{2} + 1$
split pairs. Similarly, the maximum directly below $T$ (a point of
$K_-$) is involved in $\frac{m_+}{2} + 1$ split pairs. Each of the
$\frac{m_++m_-}{2}+2$ split pairs reduces the width by 4, giving a
total of $2m+8$.

\item  As in Case (1) the minimum directly above $T$ and {\tt min} are
involved in $\frac{m_-}{2}$ split pairs each. As in Case (2) the
maximum directly below $T$ (a point of $K_-$) is involved in
$\frac{m_+}{2} + 1$ split pairs. Together we get $m_- +
\frac{m_+}{2} + 1$ split pairs, yielding a reduction in width of
at least $4 m_- + 2 m_+ + 4 \geq 2 m + 4$.

\item Symmetric to case (3) we get a reduction of at least $2m+4$.
\end{enumerate}

Next we reattach the lowest minimum of $K_+$ to the highest
maximum in $K_-$ to obtain a presentation of $K$.  This will
reduce width by exactly two. See Figure \ref{fig:Thinner}.

Splitting to form the link increased the width by $2m+2$ and the
final reattachment reduced it by $2$, a net increase of $2m$. In
cases 2, 3, and 4 the separation of $K_+$ and $K_-$ reduced the
width by at least $2m+8$ (case 2) or $2m+4$ (cases 3 and 4). In
these cases, we obtain an overall reduction of at least 4,
contradicting thin position. In Case 1 we lowered the width by
$4m$ which yields an overall reduction unless $m=0$, our desired
conclusion. (Note that $m=0$ means that $K_+$ and $K_-$ do not
cross the level $T$ and hence force us to be in case 1.)
\end{proof}

This completes the proof of Theorem \ref{thm:connect-small-knots}.
\end{proof}

\bibliographystyle{gtart}

\Addresses\recd

\end{document}